\documentclass[12pt,thmsa]{article}%
\usepackage{amssymb}
\usepackage{amsfonts}
\usepackage{sw20bams}%
\usepackage{amsmath}%
\usepackage{graphicx}
\providecommand{\U}[1]{\protect\rule{.1in}{.1in}}
\begin{document}

\author{Steven R. Finch}
\title{Mean Width of a Regular Cross-Polytope}
\date{March 12, 2016}
\maketitle

\begin{abstract}
The expected range of a sample of $n+1$ normally distributed variables is
known to be related to the mean width of a regular $n$-simplex. We show that
the expected maximum $\mu_{n}$ of a sample of $n$ half-normally distributed
variables is related to the mean width of a regular $n$-crosspolytope. Both of
these relations have mean square counterparts. An expression for $\mu_{5} $ is
found and is believed to be new.

\end{abstract}

\footnotetext{Copyright \copyright \ 2011, 2016 by Steven R. Finch. All rights
reserved.}Let $C$ be a convex body in $\mathbb{R}^{n}$. A\ \textbf{width} is
the distance between a pair of parallel $C$-supporting planes (linear
varieties of dimension $n-1$). Every unit vector $u\in\mathbb{R}^{n}$
determines a unique such pair of planes orthogonal to $u$ and hence a width
$w(u)$. Let $u $ be uniformly distributed on the unit sphere $S^{n-1}\subset$
$\mathbb{R}^{n}$. Then $w$ is a random variable and
\[
\mathbb{E}\left(  w_{3}\right)  =\frac{3}{\pi}\arccos\left(  \frac{1}%
{3}\right)
\]
for $C=$ the regular $3$-crosspolytope (octahedron) in $\mathbb{R}^{3}$ with
edges of unit length. \ Our contribution is to extend the preceding
\textbf{mean width} result to regular $n$-crosspolytopes in $\mathbb{R}^{n}$
for $n\leq6$, starting with an integral formula in \cite{BH, HZ}. We similarly
extend the following \textbf{mean square width} result:%
\[
\mathbb{E}\left(  w_{3}^{2}\right)  =\frac{2}{3}\left(  1+\frac{2\sqrt{3}}%
{\pi}\right)
\]
which, as far as is known, first appeared in \cite{Fi1}.

\section{Order Statistics}

Let $X_{1},$ $X_{2},$ $...,$ $X_{n}$ denote a random sample obtained by
\textquotedblleft folding\textquotedblright\ Normal $(0,1)$ across the origin,
that is, with density function $f$ and cumulative distribution $F$:%
\[%
\begin{array}
[c]{ccc}%
f(x)=\sqrt{\dfrac{2}{\pi}}\exp\left(  -\dfrac{x^{2}}{2}\right)  , &  & F(x)=%
%TCIMACRO{\dint \limits_{0}^{x}}%
%BeginExpansion
{\displaystyle\int\limits_{0}^{x}}
%EndExpansion
f(\xi)\,d\xi=\operatorname*{erf}\left(  \dfrac{x}{\sqrt{2}}\right)
\end{array}
\]
for $x>0$ only.

The first two moments of the half-normal \textbf{maximum}
\[
m_{n}=\max\{X_{1},X_{2},...,X_{n}\}
\]
are given by \cite{Da}%
\[%
\begin{array}
[c]{ccc}%
\mu_{n}=\mathbb{E}(m_{n})=n%
%TCIMACRO{\dint \limits_{0}^{\infty}}%
%BeginExpansion
{\displaystyle\int\limits_{0}^{\infty}}
%EndExpansion
x\,F(x)^{n-1}f(x)dx, &  & \nu_{n}=\mathbb{E}(m_{n}^{2})=n%
%TCIMACRO{\dint \limits_{0}^{\infty}}%
%BeginExpansion
{\displaystyle\int\limits_{0}^{\infty}}
%EndExpansion
x^{2}\,F(x)^{n-1}f(x)dx.
\end{array}
\]
For small $n$, exact expressions are possible \cite{Gv1, Kz}:
\[%
\begin{array}
[c]{lll}%
\mu_{2}=\frac{2}{\sqrt{\pi}}=1.128..., &  &
\begin{array}
[c]{c}%
\nu_{2}=1+\frac{2}{\pi}=1.636...,
\end{array}
\\
\mu_{3}=\frac{12}{\sqrt{\pi}}S_{2}=1.326..., &  &
\begin{array}
[c]{c}%
\nu_{3}=1+\frac{2\sqrt{3}}{\pi}=2.102...,
\end{array}
\\
\mu_{4}=\frac{12}{\sqrt{\pi}}\left(  1-4S_{2}\right)  =1.464..., &  &
\begin{array}
[c]{c}%
\nu_{4}=1+\frac{8\sqrt{3}}{3\pi}=2.470...,
\end{array}
\\
\mu_{5}=\frac{5}{\sqrt{\pi}}\left(  -\sqrt{2}+8S_{2}+16\sqrt{2}T_{1}^{\prime
}\right)  =1.569..., &  &
\begin{array}
[c]{c}%
\nu_{5}=1+\frac{20\sqrt{3}}{\pi}\left(  1-4S_{3}\right)  =2.773...,
\end{array}
\\
\mu_{6}=\frac{30}{\sqrt{\pi}}\left(  1-8S_{2}+16T_{2}\right)  =1.653... &  &
\begin{array}
[c]{c}
\end{array}
\end{array}
\]
where \cite{Wa}
\[
S_{k}=\frac{\sqrt{k}}{\pi}%
%TCIMACRO{\dint \limits_{0}^{\pi/4}}%
%BeginExpansion
{\displaystyle\int\limits_{0}^{\pi/4}}
%EndExpansion
\frac{dx}{\sqrt{k+\sec(x)^{2}}}=\frac{1}{2\pi}\operatorname{arcsec}\left(
k+1\right)  ,
\]%
\[
T_{2}=\frac{\sqrt{2}}{\pi^{2}}%
%TCIMACRO{\dint \limits_{0}^{\pi/4}}%
%BeginExpansion
{\displaystyle\int\limits_{0}^{\pi/4}}
%EndExpansion
\,%
%TCIMACRO{\dint \limits_{0}^{\pi/4}}%
%BeginExpansion
{\displaystyle\int\limits_{0}^{\pi/4}}
%EndExpansion
\frac{dx\,dy}{\sqrt{2+\sec(x)^{2}+\sec(y)^{2}}}=\frac{1}{2\pi^{2}}%
%TCIMACRO{\dint \limits_{0}^{\pi\,S_{2}}}%
%BeginExpansion
{\displaystyle\int\limits_{0}^{\pi\,S_{2}}}
%EndExpansion
\operatorname{arcsec}\left(  1+\frac{6}{2-\tan(z)^{2}}\right)  dz,
\]%
\begin{align*}
T_{1}^{\prime}  & =\frac{1}{\pi^{2}}%
%TCIMACRO{\dint \limits_{0}^{\pi/4}}%
%BeginExpansion
{\displaystyle\int\limits_{0}^{\pi/4}}
%EndExpansion
\,%
%TCIMACRO{\dint \limits_{0}^{\pi/4}}%
%BeginExpansion
{\displaystyle\int\limits_{0}^{\pi/4}}
%EndExpansion
\frac{dx\,dy}{1+\sec(x)^{2}+\sec(y)^{2}}\\
& =\frac{1-2\sqrt{2}S_{2}}{16}+\frac{1}{2\sqrt{2}\pi^{2}}%
%TCIMACRO{\dint \limits_{0}^{\pi\,S_{2}}}%
%BeginExpansion
{\displaystyle\int\limits_{0}^{\pi\,S_{2}}}
%EndExpansion
\,\,\frac{\operatorname{arcsec}\left(  2-\frac{3}{1+4\cot(z)^{2}}\right)
}{1-\frac{6}{2-\tan(z)^{2}}}dz.
\end{align*}
The preceding table complements an analogous table in \cite{Fi2} for first and
second moments of the normal maximum (rather than half-normal). Similar
expressions for $\nu_{6}=3.032...$ and $\mu_{7}=1.723...$ remain to be found.
\ Our expression for $\mu_{5}$ is believed to be new and details will be
provided shortly.

\section{Starting Point}

A\ regular $n$-crosspolytope with edges of unit length has mean width
\cite{BH, HZ}
\[
\mathbb{E}\left(  w_{n}\right)  =\frac{2\sqrt{2}n(n-1)}{\pi}\frac
{\Gamma\left(  \frac{n}{2}\right)  }{\Gamma\left(  \frac{n+1}{2}\right)  }%
%TCIMACRO{\dint \limits_{0}^{\infty}}%
%BeginExpansion
{\displaystyle\int\limits_{0}^{\infty}}
%EndExpansion
e^{-2x^{2}}\operatorname{erf}(x)^{n-2}\,dx.
\]
Integrating by parts yields%
\begin{align*}%
%TCIMACRO{\dint \limits_{0}^{\infty}}%
%BeginExpansion
{\displaystyle\int\limits_{0}^{\infty}}
%EndExpansion
e^{-2x^{2}}\operatorname{erf}(x)^{n-2}\,dx  & =\frac{\sqrt{\pi}}{n-1}%
%TCIMACRO{\dint \limits_{0}^{\infty}}%
%BeginExpansion
{\displaystyle\int\limits_{0}^{\infty}}
%EndExpansion
x\,e^{-x^{2}}\operatorname{erf}(x)^{n-1}\,dx\\
& =\frac{\sqrt{\pi}}{2(n-1)}%
%TCIMACRO{\dint \limits_{0}^{\infty}}%
%BeginExpansion
{\displaystyle\int\limits_{0}^{\infty}}
%EndExpansion
x\,e^{-x^{2}/2}\operatorname{erf}\left(  \frac{x}{\sqrt{2}}\right)
^{n-1}\,dx\\
& =\frac{\pi}{2\sqrt{2}(n-1)}%
%TCIMACRO{\dint \limits_{0}^{\infty}}%
%BeginExpansion
{\displaystyle\int\limits_{0}^{\infty}}
%EndExpansion
x\,f(x)F(x)^{n-1}\,dx.
\end{align*}
We recognize the latter integral as $\mu_{n}/n$; hence
\[
\mathbb{E}\left(  w_{n}\right)  =\frac{\Gamma\left(  \frac{n}{2}\right)
}{\Gamma\left(  \frac{n+1}{2}\right)  }\mu_{n}
\]
and conjecturally
\[
\mathbb{E}\left(  w_{n}^{2}\right)  =\frac{2}{n}\nu_{n}.
\]
Numerical confirmation for $n\leq6$ is possible via the computer algebra
technique described in \cite{Fi1}.

In summary, we have mean width results%
\[
\mathbb{E}\left(  w_{2}\right)  =\frac{4}{\pi}=1.273239544735162...,
\]%
\[
\mathbb{E}\left(  w_{3}\right)  =6S_{2}=1.175479656091821...,
\]%
\[
\mathbb{E}\left(  w_{4}\right)  =\frac{16}{\pi}\left(  1-4S_{2}\right)
=1.101845693159859...,
\]%
\[
\mathbb{E}\left(  w_{5}\right)  =\frac{15}{8}\left(  -\sqrt{2}+8S_{2}%
+16\sqrt{2}T_{1}^{\prime}\right)  =1.043421681509785...,
\]%
\[
\mathbb{E}\left(  w_{6}\right)  =\frac{32}{\pi}\left(  1-8S_{2}+16T_{2}%
\right)  =0.995378656038812...
\]
and mean square width results%
\[
\mathbb{E}\left(  w_{2}^{2}\right)  =1+\frac{2}{\pi}=1.636619772367581...,
\]%
\[
\mathbb{E}\left(  w_{3}^{2}\right)  =\frac{2}{3}\left(  1+\frac{2\sqrt{3}}%
{\pi}\right)  =1.401771860562389...,
\]%
\[
\mathbb{E}\left(  w_{4}^{2}\right)  =\frac{1}{2}\left(  1+\frac{8\sqrt{3}%
}{3\pi}\right)  =1.235105193895722...,
\]%
\[
\mathbb{E}\left(  w_{5}^{2}\right)  =\frac{2}{5}\left(  1+\frac{20\sqrt{3}%
}{\pi}\left(  1-4S_{3}\right)  \right)  =1.109499626837713....
\]

\section{Details Underlying $\mu_{6}$}

To compute $\mu_{6}$ is less difficult than to compute $\mu_{5}$. \ One
approach is to use \cite{Gv2, BG, AL}, relating moments of order statistics
for a symmetric distribution with moments of the corresponding folded
distribution. \ We will follow \cite{Wa}, however, because no other approach
seems to be applicable to evaluating $\mu_{5}$ as well.

As a prelude,%
\[
t\,f(t)=-\frac{d}{dt}f(t)
\]
hence%
\[
\mu_{6}=6%
%TCIMACRO{\dint \limits_{0}^{\infty}}%
%BeginExpansion
{\displaystyle\int\limits_{0}^{\infty}}
%EndExpansion
t\,F(t)^{5}f(t)dt=-6%
%TCIMACRO{\dint \limits_{0}^{\infty}}%
%BeginExpansion
{\displaystyle\int\limits_{0}^{\infty}}
%EndExpansion
F(t)^{5}\frac{d}{dt}f(t)dt=30%
%TCIMACRO{\dint \limits_{0}^{\infty}}%
%BeginExpansion
{\displaystyle\int\limits_{0}^{\infty}}
%EndExpansion
F(t)^{4}f(t)^{2}dt
\]
via integration by parts. Also
\begin{align*}
F(t)^{2}  & =\dfrac{2}{\pi}%
%TCIMACRO{\dint \limits_{0}^{t}}%
%BeginExpansion
{\displaystyle\int\limits_{0}^{t}}
%EndExpansion%
%TCIMACRO{\dint \limits_{0}^{t}}%
%BeginExpansion
{\displaystyle\int\limits_{0}^{t}}
%EndExpansion
\exp\left(  -\dfrac{x^{2}+y^{2}}{2}\right)  dx\,dy=\dfrac{4}{\pi}%
%TCIMACRO{\dint \limits_{0}^{\pi/4}}%
%BeginExpansion
{\displaystyle\int\limits_{0}^{\pi/4}}
%EndExpansion
\,\,\,%
%TCIMACRO{\dint \limits_{0}^{t\sec(\theta)}}%
%BeginExpansion
{\displaystyle\int\limits_{0}^{t\sec(\theta)}}
%EndExpansion
\exp\left(  -\dfrac{r^{2}}{2}\right)  r\,dr\,d\theta\\
& =1-\dfrac{4}{\pi}%
%TCIMACRO{\dint \limits_{0}^{\pi/4}}%
%BeginExpansion
{\displaystyle\int\limits_{0}^{\pi/4}}
%EndExpansion
\exp\left(  -\dfrac{t^{2}\sec(\theta)^{2}}{2}\right)  d\theta
\end{align*}
via transformation to polar coordinates.

We deduce that \
\begin{align*}
\frac{\mu_{6}}{30}  & =%
%TCIMACRO{\dint \limits_{0}^{\infty}}%
%BeginExpansion
{\displaystyle\int\limits_{0}^{\infty}}
%EndExpansion
\left[  1-\dfrac{4}{\pi}%
%TCIMACRO{\dint \limits_{0}^{\pi/4}}%
%BeginExpansion
{\displaystyle\int\limits_{0}^{\pi/4}}
%EndExpansion
\exp\left(  -\dfrac{t^{2}\sec(\theta)^{2}}{2}\right)  d\theta\right]
^{2}f(t)^{2}dt\\
& =%
%TCIMACRO{\dint \limits_{0}^{\infty}}%
%BeginExpansion
{\displaystyle\int\limits_{0}^{\infty}}
%EndExpansion
f(t)^{2}dt-\frac{8}{\pi}%
%TCIMACRO{\dint \limits_{0}^{\infty}}%
%BeginExpansion
{\displaystyle\int\limits_{0}^{\infty}}
%EndExpansion
\,%
%TCIMACRO{\dint \limits_{0}^{\pi/4}}%
%BeginExpansion
{\displaystyle\int\limits_{0}^{\pi/4}}
%EndExpansion
\exp\left(  -\dfrac{t^{2}\sec(\theta)^{2}}{2}\right)  f(t)^{2}d\theta\,dt\\
& \;\;\;\;\;\;\;\;\;+\frac{16}{\pi^{2}}%
%TCIMACRO{\dint \limits_{0}^{\infty}}%
%BeginExpansion
{\displaystyle\int\limits_{0}^{\infty}}
%EndExpansion
\,%
%TCIMACRO{\dint \limits_{0}^{\pi/4}}%
%BeginExpansion
{\displaystyle\int\limits_{0}^{\pi/4}}
%EndExpansion
\,\,%
%TCIMACRO{\dint \limits_{0}^{\pi/4}}%
%BeginExpansion
{\displaystyle\int\limits_{0}^{\pi/4}}
%EndExpansion
\exp\left(  -\dfrac{t^{2}\left(  \sec(\varphi)^{2}+\sec(\psi)^{2}\right)  }%
{2}\right)  f(t)^{2}d\varphi\,d\psi\,dt\\
& =\frac{1}{\sqrt{\pi}}-\frac{8}{\sqrt{\pi}}S_{2}+\frac{32}{\pi^{3}}%
%TCIMACRO{\dint \limits_{0}^{\infty}}%
%BeginExpansion
{\displaystyle\int\limits_{0}^{\infty}}
%EndExpansion
\,%
%TCIMACRO{\dint \limits_{0}^{\pi/4}}%
%BeginExpansion
{\displaystyle\int\limits_{0}^{\pi/4}}
%EndExpansion
\,\,%
%TCIMACRO{\dint \limits_{0}^{\pi/4}}%
%BeginExpansion
{\displaystyle\int\limits_{0}^{\pi/4}}
%EndExpansion
\exp\left(  -\dfrac{t^{2}\left(  2+\sec(\varphi)^{2}+\sec(\psi)^{2}\right)
}{2}\right)  d\varphi\,d\psi\,dt\\
& =\frac{1}{\sqrt{\pi}}-\frac{8}{\sqrt{\pi}}S_{2}+\frac{32}{\pi^{3}}%
\sqrt{\frac{\pi}{2}}%
%TCIMACRO{\dint \limits_{0}^{\pi/4}}%
%BeginExpansion
{\displaystyle\int\limits_{0}^{\pi/4}}
%EndExpansion
\,\,%
%TCIMACRO{\dint \limits_{0}^{\pi/4}}%
%BeginExpansion
{\displaystyle\int\limits_{0}^{\pi/4}}
%EndExpansion
\frac{d\varphi\,d\psi}{\sqrt{2+\sec(\varphi)^{2}+\sec(\psi)^{2}}}\\
& =\frac{1}{\sqrt{\pi}}-\frac{8}{\sqrt{\pi}}S_{2}+\frac{32}{\pi^{3}}%
\sqrt{\frac{\pi}{2}}\frac{\pi^{2}}{\sqrt{2}}T_{2}=\frac{1}{\sqrt{\pi}}\left(
1-8S_{2}+16T_{2}\right)  .
\end{align*}
As an interlude,%
\[%
%TCIMACRO{\dint \limits_{0}^{\pi/4}}%
%BeginExpansion
{\displaystyle\int\limits_{0}^{\pi/4}}
%EndExpansion
\frac{d\psi}{\sqrt{\kappa+\sec(\psi)^{2}}}=\frac{\pi}{\sqrt{\kappa}}S_{\kappa
}=\frac{1}{2\sqrt{\kappa}}\operatorname{arcsec}\left(  \kappa+1\right)
\]
hence%
\[%
%TCIMACRO{\dint \limits_{0}^{\pi/4}}%
%BeginExpansion
{\displaystyle\int\limits_{0}^{\pi/4}}
%EndExpansion
\,\,%
%TCIMACRO{\dint \limits_{0}^{\pi/4}}%
%BeginExpansion
{\displaystyle\int\limits_{0}^{\pi/4}}
%EndExpansion
\frac{d\varphi\,d\psi}{\sqrt{2+\sec(\varphi)^{2}+\sec(\psi)^{2}}}=\frac{1}{2}%
%TCIMACRO{\dint \limits_{0}^{\pi/4}}%
%BeginExpansion
{\displaystyle\int\limits_{0}^{\pi/4}}
%EndExpansion
\frac{1}{\sqrt{2+\sec(\varphi)^{2}}}\operatorname{arcsec}\left(
3+\sec(\varphi)^{2}\right)  d\varphi
\]
upon substituting $2+\sec(\varphi)^{2}$ for $\kappa$.\ Now let
\[
2\sin(\varphi)^{2}=3\sin(\theta)^{2},
\]
then%
\[
2-2\cos(\varphi)^{2}=3-3\cos(\theta)^{2}
\]
which implies that%
\[%
\begin{array}
[c]{ccc}%
\sin(\theta)^{2}=\dfrac{2}{3}\sin(\varphi)^{2}, &  & \cos(\theta)^{2}%
=\dfrac{1}{3}\left(  1+2\cos(\varphi)^{2}\right)
\end{array}
\]
hence%
\begin{align*}
4\sin(\varphi)\cos(\varphi)d\varphi & =6\sin(\theta)\cos(\theta)d\theta\\
& =\frac{6}{3}\sqrt{\left(  2\sin(\varphi)^{2}\right)  \left(  1+2\cos
(\varphi)^{2}\right)  }\,d\theta\\
& =2\sqrt{2}\sin(\varphi)\sqrt{1+2\cos(\varphi)^{2}}\,d\theta
\end{align*}
hence%
\[
d\varphi=\frac{1}{\sqrt{2}}\frac{\sqrt{1+2\cos(\varphi)^{2}}}{\cos(\varphi
)}d\theta=\frac{1}{\sqrt{2}}\sqrt{2+\sec(\varphi)^{2}}\,d\theta.
\]
The preceding integral becomes%
\[
\frac{1}{2\sqrt{2}}%
%TCIMACRO{\dint \limits_{0}^{\pi\,S_{2}}}%
%BeginExpansion
{\displaystyle\int\limits_{0}^{\pi\,S_{2}}}
%EndExpansion
\operatorname{arcsec}\left(  3+\sec(\varphi)^{2}\right)  d\theta=\frac
{1}{2\sqrt{2}}%
%TCIMACRO{\dint \limits_{0}^{\pi\,S_{2}}}%
%BeginExpansion
{\displaystyle\int\limits_{0}^{\pi\,S_{2}}}
%EndExpansion
\operatorname{arcsec}\left(  1+\frac{6}{2-\tan(\theta)^{2}}\right)  d\theta
\]
since%
\begin{align*}
2+\sec(\varphi)^{2}  & =2+\frac{1}{1-\sin(\varphi)^{2}}=\frac{3-2\sin
(\varphi)^{2}}{1-\sin(\varphi)^{2}}=\frac{3-3\sin(\theta)^{2}}{1-\frac{3}%
{2}\sin(\theta)^{2}}=\frac{6\cos(\theta)^{2}}{2-3\sin(\theta)^{2}}\\
& =\frac{6\cos(\theta)^{2}}{2\cos(\theta)^{2}-\sin(\theta)^{2}}=\frac
{6}{2-\tan(\theta)^{2}}.
\end{align*}
Finally, as a postlude,%
\[
2\arccos(\zeta)=\arccos\left(  2\zeta^{2}-1\right)
\]
and $\varphi=\pi/4$ implies that $2=3\cos(\theta)^{2}$, thus the upper
integration limit is
\[
\theta=\arccos\left(  \sqrt{\tfrac{2}{3}}\right)  =\tfrac{1}{2}\arccos\left(
\tfrac{4}{3}-1\right)  =\tfrac{1}{2}\operatorname{arcsec}\left(  3\right)
=\pi\,S_{2}.
\]

\section{Details Underlying $\mu_{5}$}

Integration by parts is not employed this time. \ Beginning with%

\[
\mu_{5}=5%
%TCIMACRO{\dint \limits_{0}^{\infty}}%
%BeginExpansion
{\displaystyle\int\limits_{0}^{\infty}}
%EndExpansion
t\,F(t)^{4}f(t)dt
\]
it follows that \
\begin{align*}
\frac{\mu_{5}}{5}  & =%
%TCIMACRO{\dint \limits_{0}^{\infty}}%
%BeginExpansion
{\displaystyle\int\limits_{0}^{\infty}}
%EndExpansion
\left[  1-\dfrac{4}{\pi}%
%TCIMACRO{\dint \limits_{0}^{\pi/4}}%
%BeginExpansion
{\displaystyle\int\limits_{0}^{\pi/4}}
%EndExpansion
\exp\left(  -\dfrac{t^{2}\sec(\theta)^{2}}{2}\right)  d\theta\right]
^{2}t\,f(t)dt\\
& =%
%TCIMACRO{\dint \limits_{0}^{\infty}}%
%BeginExpansion
{\displaystyle\int\limits_{0}^{\infty}}
%EndExpansion
t\,f(t)dt-\frac{8}{\pi}%
%TCIMACRO{\dint \limits_{0}^{\infty}}%
%BeginExpansion
{\displaystyle\int\limits_{0}^{\infty}}
%EndExpansion
\,%
%TCIMACRO{\dint \limits_{0}^{\pi/4}}%
%BeginExpansion
{\displaystyle\int\limits_{0}^{\pi/4}}
%EndExpansion
\exp\left(  -\dfrac{t^{2}\sec(\theta)^{2}}{2}\right)  t\,f(t)d\theta\,dt\\
& \;\;\;\;\;\;\;\;\;+\frac{16}{\pi^{2}}%
%TCIMACRO{\dint \limits_{0}^{\infty}}%
%BeginExpansion
{\displaystyle\int\limits_{0}^{\infty}}
%EndExpansion
\,%
%TCIMACRO{\dint \limits_{0}^{\pi/4}}%
%BeginExpansion
{\displaystyle\int\limits_{0}^{\pi/4}}
%EndExpansion
\,\,%
%TCIMACRO{\dint \limits_{0}^{\pi/4}}%
%BeginExpansion
{\displaystyle\int\limits_{0}^{\pi/4}}
%EndExpansion
\exp\left(  -\dfrac{t^{2}\left(  \sec(\varphi)^{2}+\sec(\psi)^{2}\right)  }%
{2}\right)  t\,f(t)d\varphi\,d\psi\,dt\\
& =\sqrt{\dfrac{2}{\pi}}-\frac{2\sqrt{2}-8S_{2}}{\sqrt{\pi}}+\frac{16}{\pi
^{2}}\sqrt{\dfrac{2}{\pi}}%
%TCIMACRO{\dint \limits_{0}^{\infty}}%
%BeginExpansion
{\displaystyle\int\limits_{0}^{\infty}}
%EndExpansion
\,%
%TCIMACRO{\dint \limits_{0}^{\pi/4}}%
%BeginExpansion
{\displaystyle\int\limits_{0}^{\pi/4}}
%EndExpansion
\,\,%
%TCIMACRO{\dint \limits_{0}^{\pi/4}}%
%BeginExpansion
{\displaystyle\int\limits_{0}^{\pi/4}}
%EndExpansion
t\exp\left(  -\dfrac{t^{2}\left(  1+\sec(\varphi)^{2}+\sec(\psi)^{2}\right)
}{2}\right)  d\varphi\,d\psi\,dt\\
& =\frac{1}{\sqrt{\pi}}\left\{  \sqrt{2}-(2\sqrt{2}-8S_{2})+\frac{16\sqrt{2}%
}{\pi^{2}}%
%TCIMACRO{\dint \limits_{0}^{\pi/4}}%
%BeginExpansion
{\displaystyle\int\limits_{0}^{\pi/4}}
%EndExpansion
\,\,%
%TCIMACRO{\dint \limits_{0}^{\pi/4}}%
%BeginExpansion
{\displaystyle\int\limits_{0}^{\pi/4}}
%EndExpansion
\frac{d\varphi\,d\psi}{1+\sec(\varphi)^{2}+\sec(\psi)^{2}}\right\} \\
& =\frac{1}{\sqrt{\pi}}\left(  -\sqrt{2}+8S_{2}+16\sqrt{2}T_{1}^{\prime
}\right)  .
\end{align*}
As an interlude,%
\[%
%TCIMACRO{\dint \limits_{0}^{\pi/4}}%
%BeginExpansion
{\displaystyle\int\limits_{0}^{\pi/4}}
%EndExpansion
\frac{d\psi}{\kappa+\sec(\psi)^{2}}=\frac{\pi}{4\kappa}-\frac
{\operatorname{arcsec}\left(  1+\frac{2}{\kappa}\right)  }{2\kappa\sqrt
{\kappa+1}}
\]
hence%
\begin{align*}%
%TCIMACRO{\dint \limits_{0}^{\pi/4}}%
%BeginExpansion
{\displaystyle\int\limits_{0}^{\pi/4}}
%EndExpansion
\,\,%
%TCIMACRO{\dint \limits_{0}^{\pi/4}}%
%BeginExpansion
{\displaystyle\int\limits_{0}^{\pi/4}}
%EndExpansion
\frac{d\varphi\,d\psi}{1+\sec(\varphi)^{2}+\sec(\psi)^{2}}  & =%
%TCIMACRO{\dint \limits_{0}^{\pi/4}}%
%BeginExpansion
{\displaystyle\int\limits_{0}^{\pi/4}}
%EndExpansion
\left[  \frac{\pi}{4\left(  1+\sec(\varphi)^{2}\right)  }-\frac
{\operatorname{arcsec}\left(  1+\frac{2}{1+\sec(\varphi)^{2}}\right)
}{2\left(  1+\sec(\varphi)^{2}\right)  \sqrt{2+\sec(\varphi)^{2}}}\right]
d\varphi\\
& =\frac{\pi^{2}}{16}(1-2\sqrt{2}S_{2})-%
%TCIMACRO{\dint \limits_{0}^{\pi/4}}%
%BeginExpansion
{\displaystyle\int\limits_{0}^{\pi/4}}
%EndExpansion
\frac{\operatorname{arcsec}\left(  1+\frac{2}{1+\sec(\varphi)^{2}}\right)
}{2\left(  1+\sec(\varphi)^{2}\right)  \sqrt{2+\sec(\varphi)^{2}}}d\varphi
\end{align*}
upon substituting $1+\sec(\varphi)^{2}$ for $\kappa$.\ Now, as before, let
\[
2\sin(\varphi)^{2}=3\sin(\theta)^{2}
\]
then the preceding integral becomes%
\[
\frac{1}{2\sqrt{2}}%
%TCIMACRO{\dint \limits_{0}^{\pi\,S_{2}}}%
%BeginExpansion
{\displaystyle\int\limits_{0}^{\pi\,S_{2}}}
%EndExpansion
\,\,\frac{\operatorname{arcsec}\left(  1+\frac{2}{-1+\frac{6}{2-\tan
(\theta)^{2}}}\right)  }{-1+\frac{6}{2-\tan(\theta)^{2}}}d\theta=-\frac
{1}{2\sqrt{2}}%
%TCIMACRO{\dint \limits_{0}^{\pi\,S_{2}}}%
%BeginExpansion
{\displaystyle\int\limits_{0}^{\pi\,S_{2}}}
%EndExpansion
\,\,\frac{\operatorname{arcsec}\left(  2-\frac{3}{1+4\cot(\theta)^{2}}\right)
}{1-\frac{6}{2-\tan(\theta)^{2}}}d\theta
\]
since%
\begin{align*}
-1+\frac{2}{-1+\frac{6}{2-\tan(\theta)^{2}}}  & =\frac{-\left(  -1+\frac
{6}{2-\tan(\theta)^{2}}\right)  +2}{-1+\frac{6}{2-\tan(\theta)^{2}}}%
=\frac{3-\frac{6}{2-\tan(\theta)^{2}}}{-1+\frac{6}{2-\tan(\theta)^{2}}}%
=\frac{6-3\tan(\theta)^{2}-6}{-2+\tan(\theta)^{2}+6}\\
& =\frac{-3\tan(\theta)^{2}}{4+\tan(\theta)^{2}}=-\frac{3}{1+4\cot(\theta
)^{2}}.
\end{align*}
Everything comes together at this point, completing the derivation.

\section{Asymptotics}

We turn now to the asymptotic distribution of $m_{n}$ as $n\rightarrow\infty$.
Define $a_{n}$ to be the positive solution of the equation \cite{Gb, Ha}%
\[
2\pi\,a_{n}^{2}\exp\left(  a_{n}^{2}\right)  =n^{2},
\]
that is,%
\[
a_{n}=\sqrt{W\left(  \frac{n^{2}}{2\pi}\right)  }\sim\sqrt{2\ln(n)}-\frac
{1}{2}\frac{\ln(\ln(n))+\ln(4\pi)}{\sqrt{2\ln(n)}}
\]
in terms of the Lambert $W$ function \cite{Cr}. \ An adjustment is needed that
corresponds to folding (truncating). Define \cite{Gb}%
\[
a_{n}^{\prime}=a_{n}+\frac{1}{a_{n}}\ln(2)\sim\sqrt{2\ln(n)}-\frac{1}{2}%
\frac{\ln(\ln(n))+\ln(\pi)}{\sqrt{2\ln(n)}},
\]
then the required density is
\[
\lim_{n\rightarrow\infty}\frac{d}{dy}\operatorname*{P}\left(  \sqrt{2\ln
(n)}(m_{n}-a_{n}^{\prime})<y\right)  =\exp\left(  -e^{-y}\right)  .
\]
A random variable $Y$, distributed as such, satisfies
\[%
\begin{array}
[c]{ccc}%
\mathbb{E}(Y)=\gamma, &  & \mathbb{E}(Y^{2})=\dfrac{\pi^{2}}{6}+\gamma^{2}%
\end{array}
\]
where $\gamma$ is the Euler-Mascheroni constant \cite{Fi3}. \ This implies
that
\[
\mu_{n}\sim a_{n}^{\prime}+\frac{\gamma}{\sqrt{2\ln(n)}}\sim\sqrt{2\ln
(n)}-\frac{1}{2}\frac{\ln(\ln(n))+\ln(\pi)-2\gamma}{\sqrt{2\ln(n)}}
\]
and hence%
\begin{align*}
\mathbb{E}(w_{n})  & =\frac{\Gamma\left(  \frac{n}{2}\right)  }{\Gamma\left(
\frac{n+1}{2}\right)  }\cdot\mu_{n}\\
& \sim\frac{1}{\sqrt{2n}}\left(  1+\frac{1}{4n}\right)  \cdot\left(
a_{n}^{\prime}+\frac{\gamma}{\sqrt{2\ln(n)}}\right) \\
& \sim\sqrt{\frac{\ln(n)}{n}}-\frac{\ln(\ln(n))+\ln(\pi)-2\gamma}{4\sqrt
{n\ln(n)}}.
\end{align*}
More terms in the asymptotic expansion are possible. The first term of this
approximation is consistent with \cite{Vs}.

If we rescale length so that the inradius of the $n$-crosspolytope is $1$ and
denote adjusted width by $\widehat{w}_{n}$, then%
\[
\mathbb{E}(\widehat{w}_{n})\sim\sqrt{2n}\cdot\sqrt{\frac{\ln(n)}{n}}\sim
\sqrt{2\ln(n)}
\]
because, in our original scaling, the inradius is $\sqrt{1/(2n)}$. \ 

After completing this work, we discovered an article \cite{Bl} on metric
properties of the only three regular polytopes in $5$-dimensional space or
higher. \ Our paper and its predecessor \cite{Fi1} serve, in some sense, to
complement this article.

\bigskip
\end{document}